\documentclass[reqno,12pt]{amsart}
\usepackage{a4wide}
\usepackage{amscd}
\usepackage{amsmath}
\usepackage{amsfonts}
\usepackage{amssymb}
\usepackage{multirow}
\usepackage{color}
\input xy
\xyoption{all}

\usepackage[OT2,OT1]{fontenc}
\newcommand\cyr{%
\renewcommand\rmdefault{wncyr}%
\renewcommand\sfdefault{wncyss}%
\renewcommand\encodingdefault{OT2}%
\normalfont
\selectfont}
\DeclareTextFontCommand{\textcyr}{\cyr}

\newcommand{\ce}{{\mathfrak{e}}}

\newcommand{\cf}{{\mathfrak{f}}}

\newcommand{\cg}{{\mathfrak{g}}}

\newcommand{\ch}{{\mathfrak{h}}}

\newcommand{\co}{{\mathfrak{o}}}

\newcommand{\cp}{{\mathfrak{p}}}

\newcommand{\cs}{{\mathfrak{s}}}

\newcommand{\cu}{{\mathfrak{u}}}

\newcommand{\RR}{\mathbb{R}}

\newcommand{\Aut}{\mbox{Aut}}
\newcommand{\codim}{\mbox{codim}}
\newcommand{\Fix}{\mbox{Fix}}
\newcommand{\rk}{\mbox{rk}}

\newtheorem{thm}{Theorem}[section]

\newtheorem{cor}[thm]{Corollary}

\numberwithin{equation}{section}

\begin {document}

\title{The index of compact simple Lie groups}

\author{J\"{u}rgen Berndt}
\address{King's College London, Department of Mathematics, London WC2R 2LS, United Kingdom}
%\curraddr{}
\email{jurgen.berndt@kcl.ac.uk}
\thanks{}

\author{Carlos Olmos}
\address{Facultad de Matem\'atica, Astronom\'ia y F\'isica, Universidad Nacional de C\'ordoba, 
Ciudad Universitaria, 5000 C\'ordoba, Argentina}
%\curraddr{}
\email{olmos@famaf.unc.edu.ar}
%\thanks{}

\subjclass[2010]{Primary 53C35; Secondary 53C40}

\hyphenation{geo-metry}

\begin{abstract}
Let $M$ be an irreducible Riemannian symmetric space. The index $i(M)$ of $M$ is the 
minimal codimension of a (non-trivial) totally geodesic submanifold of $M$. The purpose of this  
note is to determine the index $i(M)$ for all irreducible Riemannian symmetric spaces $M$ of type (II) 
and (IV).
\end{abstract}

\maketitle 

\section {Introduction}

Let $M$ be a connected Riemannian manifold and denote by ${\mathcal S}$ the set of all connected totally geodesic submanifolds $\Sigma$ of $M$ with $\dim(\Sigma) < \dim(M)$. The index $i(M)$ of $M$ is defined by
\[
i(M) = \min\{ \dim(M) - \dim(\Sigma) : \Sigma \in {\mathcal S}\} = \min\{ \codim(\Sigma) : \Sigma \in {\mathcal S}\}.
\]
This notion was introduced by Onishchik in \cite{On}, who also classified the irreducible simply connected Riemannian symmetric spaces $M$  with $i(M) \leq 2$. 

In \cite{BO} we investigated $i(M)$ for irreducible Riemannian symmetric spaces $M$.  We proved that the rank $\rk(M)$ of $M$ is always less than or equal to the index of $M$ and classified all irreducible Riemannian symmetric spaces $M$ with $i(M) \leq 3$. 

A totally geodesic submanifold $\Sigma$ of $M$ is called reflective if $\Sigma$ is the connected component of the fixed point set of an isometric involution on $M$. Reflective submanifolds of irreducible, simply connected Riemannian symmetric spaces of compact type were classified by Leung in \cite{L1} and \cite{L2}. Denote by ${\mathcal S}_r$ the set of all connected reflective submanifolds $\Sigma$ of $M$ with $\dim(\Sigma) < \dim(M)$. The reflective index $i_r(M)$ of $M$ is defined by
\[
i_r(M) = \min\{ \dim(M) - \dim(\Sigma) : \Sigma \in {\mathcal S}_r\} = \min\{ \codim(\Sigma) :  \Sigma \in {\mathcal S}_r\}.
\]
It is clear that $i(M) \leq i_r(M)$ and thus $i_r(M)$ is an upper bound for $i(M)$. Moreover, from \cite{L1} and \cite{L2} we can calculate $i_r(M)$ explicitly for each irreducible Riemannian symmetric space. This was done explicitly in \cite{BO2}, where we
conjectured that $i(M) = i_r(M)$ if and only if $M \neq G_2^2/SO_4$. We also verified this conjecture for a number of symmetric spaces. The purpose of this brief note is to give an affirmative answer to this conjecture for irreducible Riemannian symmetric spaces of type (II) and (IV). Since totally geodesic submanifolds are preserved under duality between symmetric spaces of compact type and of noncompact type, we can assume that $M$ is of compact type. 

Our main result is as follows. 

\begin{thm}\label{main}
Let $G$ be a simply connected, compact real simple Lie group equipped with the bi-invariant Riemannian metric induced from the Killing form of the Lie algebra $\cg$ of $G$. Then $i(G) = i_r(G)$. Moreover, if $\Sigma$ is a connected, complete, totally geodesic submanifold of $G$ with $\codim(\Sigma) = i(G)$, then the pair $(G,\Sigma)$ is as in Table \ref{Liegroup}.
\end{thm}

\begin{table}[h]
\caption{The index $i(G)$ of simply connected, compact simple Lie groups} 
\label{Liegroup} 
{\footnotesize\begin{tabular}{ | p{2cm}  p{3.5cm}  p{2cm}  p{2cm}  p{2cm}  |}
\hline \rule{0pt}{4mm}
\hspace{-1mm}$G$ & $\Sigma$ & $\dim(G)$ & $i(G)$ & Comments \\[1mm]
\hline \rule{0pt}{4mm}
\hspace{-2mm} 
$SU_2$ & $SU_2/S(U_1U_1)$ & $3$ & $1$ &  \\
$SU_3$ & $SU_3/SO_3$ & $8$ & $3$ &  \\
$SU_{r+1}$ & $S(U_rU_1)$ & $r(r+2)$ & $2r$ & $r \geq 4$ \\
$Spin_5$ & $Spin_4$, $SO_5/SO_2SO_3$ & $10$ & $4$ & \\
$Spin_{2r+1}$ & $Spin_{2r}$ & $r(2r+1)$ & $2r$ & $r \geq 3$ \\
$Sp_r$ & $Sp_{r-1}Sp_1$ & $r(2r+1)$ & $4r-4$ & $r \geq 3$\\
$Spin_{2r}$ & $Spin_{2r-1}$ & $r(2r-1)$ & $2r-1$ &  $r \geq 3$ \\
$E_6$ & $F_4$ & $78$ &  $26$ &\\
$E_7$ & $E_6U_1$ & $133$ & $54$ & \\
$E_8$ & $E_7Sp_1$ & $248$ & $112$ & \\
$F_4$ & $Spin_9$ & $52$ & $16$& \\
$G_2$ & $SU_3$, $G_2/SO_4$ & $14$ & $6$& \\[1mm]
\hline
\end{tabular}}
\end{table}

For $G \in \{SU_2,SU_3,Spin_5,G_2\}$ this was proved by Onishchik in \cite{On} and for $G = Spin_r$ with $r \geq 6$ this was proved by the authors in \cite{BO2}. Note that $Spin_6$ is isomorphic to $SU_4$. The result is new for $SU_r$ ($r \geq 4$), $Sp_r$ ($r \geq 3$), and the four exceptional Lie groups $E_6,E_7,E_8,F_4$.

As an immediate consequence of Theorem \ref{main} we get

\begin{cor}\label{conseq}
Let $G$ be a simply connected, compact simple Lie group equipped with the bi-invariant Riemannian metric induced from the Killing form of the Lie algebra $\cg$ of $G$. If $G \notin \{SU_2 ,SU_3\}$, then there exists a connected subgroup $H$ of $G$ such that the index of $G$ is equal to the codimension of $H$ in $G$.
\end{cor}

It is worthwhile to point out that $SU_3$ is a non-reflective totally geodesic submanifold of $G_2$, whereas all other totally geodesic submanifolds $\Sigma$ in Table \ref{Liegroup} are reflective submanifolds.

\section{Proof of Theorem \ref{main}}

Using \cite{L1} and \cite{L2}, we determined in \cite{BO2} the reflective index $i_r(M)$ of all irreducible Riemannian symmetrric spaces $M$ of noncompact type and the reflective submanifolds $\Sigma$ in $M$ for which $i_r(M) = \codim(\Sigma)$. Using duality between Riemannian symmetric spaces of noncompact type and of compact type, we obtain Table \ref{reflindex} for the reflective index $i_r(G)$ of all simply connected, compact simple Lie groups and the reflective submanifolds $\Sigma$ in $G$ for which $i_r(G) = \codim(\Sigma)$. 

\begin{table}[h]
\caption{The reflective index $i_r(G)$ of simply connected, compact simple Lie groups} 
\label{reflindex} 
{\footnotesize\begin{tabular}{ | p{2cm}  p{3.5cm}  p{2cm}  p{2cm}  p{2cm}  |}
\hline \rule{0pt}{4mm}
\hspace{-1mm}$G$ & $\Sigma$ & $\dim(G)$ & $i_r(G)$ & Comments \\[1mm]
\hline \rule{0pt}{4mm}
\hspace{-2mm} 
$SU_2$ & $SU_2/S(U_1U_1)$ & $3$ & $1$ &  \\
$SU_3$ & $SU_3/SO_3$ & $8$ & $3$ &  \\
$SU_{r+1}$ & $S(U_rU_1)$ & $r(r+2)$ & $2r$ & $r \geq 4$ \\
$Spin_5$ & $Spin_4$, $SO_5/SO_2SO_3$ & $10$ & $4$ & \\
$Spin_{2r+1}$ & $Spin_{2r}$ & $r(2r+1)$ & $2r$ & $r \geq 3$ \\
$Sp_r$ & $Sp_{r-1}Sp_1$ & $r(2r+1)$ & $4r-4$ & $r \geq 3$\\
$Spin_{2r}$ & $Spin_{2r-1}$ & $r(2r-1)$ & $2r-1$ &  $r \geq 3$ \\
$E_6$ & $F_4$ & $78$ &  $26$ &\\
$E_7$ & $E_6U_1$ & $133$ & $54$ & \\
$E_8$ & $E_7Sp_1$ & $248$ & $112$ & \\
$F_4$ & $Spin_9$ & $52$ & $16$& \\
$G_2$ & $G_2/SO_4$ & $14$ & $6$& \\[1mm]
\hline
\end{tabular}}
\end{table}

Note that Table \ref{reflindex} leads to Table \ref{Liegroup} when replacing $i_r(G)$ with $i(G)$ and adding $\Sigma = SU_3$ in the row for $G_2$. The two problems we thus need to solve for each $G$ are:

\begin{itemize}
\item[(1)] prove that there exists no non-reflective totally geodesic submanifold $\Sigma$ in $G$ with $\codim(\Sigma) < i_r(G)$;
\item[(2)] determine all non-reflective submanifolds $\Sigma$ in $G$ with $\codim(\Sigma) = i_r(G)$.
\end{itemize}

The following result is a crucial step towards the solution of the two problems:

\begin{thm}[Ikawa, Tasaki \cite{IT}] \label{IkawaTasaki}
A necessary and sufficient condition that a totally geodesic submanifold $\Sigma$ in a compact connected simple Lie group is maximal is that $\Sigma$ is a Cartan embedding or a maximal Lie subgroup. 
\end{thm}

The Cartan embeddings are defined as follows. Let $G/K$ be a Riemannian symmetric space of compact type and $\sigma \in \Aut(G)$ be an involutive automorphism of $G$ such that $\Fix(\sigma)^o \subset K \subset \Fix(\sigma)$, where 
\[
\Fix(\sigma) = \{g \in G : \sigma(g) = g\}
\]
and $\Fix(\sigma)^o$ is the identity component of $\Fix(\sigma)$. By definition, the automorphism $\sigma$ fixes all points in $K$ and the identity component $K^o$ of $K$ coincides with $\Fix(\sigma)^o$.

The Cartan map of $G/K$ into $G$ is the smooth map
\[
f : G/K \to G\ ,\ gK \mapsto \sigma(g)g^{-1}.
\]
The Cartan map $f$ is a covering map onto its image $\Sigma = f(G/K)$. 

Let $\theta \in \Aut(G)$ be the involutive automorphism on $G$ defined by inversion, that is,
\[
\theta : G \to G\ ,\ g \mapsto g^{-1}.
\]

We now define a third involutive automorphism $\rho \in \Aut(G)$ by $\rho = \theta \circ \sigma$. By definition, we have
\[
\rho(g) = \theta(\sigma(g)) = \sigma(g)^{-1} = \sigma(g^{-1})
\]
for all $g \in G$. Moreover, for all $g \in G$ we have
\begin{align*}
\rho(f(gK)) & = \rho(\sigma(g)g^{-1}) = \sigma((\sigma(g)g^{-1})^{-1}) = \sigma(g\sigma(g)^{-1}) = \sigma(g\sigma(g^{-1})) \\
& = \sigma(g)\sigma^2(g^{-1}) = \sigma(g)g^{-1} = f(gK).
\end{align*}
Thus the automorphism $\rho$ fixes all points in $\Sigma$. 

The automorphisms $\sigma,\theta,\rho \in \Aut(G)$ are involutive isometries of $G$, where $G$ is considered as a Riemannian symmetric space with a bi-invariant Riemannian metric. Geometrically, $\theta$ is the geodesic symmetry of $G$ at the identity $e \in G$ and its differential at $e$ is 
\[
d_e\theta : T_eG \to T_eG\ ,\ X \mapsto -X.
\]
The differential of $\sigma$ at $e$ is
\[
d_e\sigma : T_eG \to T_eG\ ,\ X \mapsto 
\begin{cases} 
X & \mbox{if } X \in T_eK, \\
-X & \mbox{if } X \in \nu_eK,
\end{cases}
\]
where $\nu_eK$ denotes the normal space of $K$ at $e$. This shows that $\sigma$ is the geodesic reflection of $G$ in the identity component $K^o$ of $K$. In particular, $K^o$ (and hence also $K$) is a totally geodesic submanifold of $G$. Since $\rho = \theta \circ \sigma$, the differential of $\rho$ at $e$ is
\[
d_e\rho : T_eG \to T_eG\ ,\ X \mapsto 
\begin{cases} 
X & \mbox{if } X \in \nu_eK, \\
-X & \mbox{if } X \in T_eK,
\end{cases}
\]
It follows that there exists a connected, complete, totally geodesic submanifold $N$ of $G$ with $e \in N$ and $T_eN = \nu_eK$. We saw above that $\rho$ fixes all points in $\Sigma$, which implies $\Sigma \subset N$ since $\Sigma$ is connected. Moreover, since $\dim(\Sigma) = \dim(G) - \dim(K) = \codim(K) = \dim(N)$ and $\Sigma$ is complete we get $\Sigma = N$. It follows that $\Sigma$ is a totally geodesic submanifold of $G$. In fact, we have proved that both $K^o$ and $\Sigma$ are reflective submanifolds of $G$ which are perpendicular to each other at $e$. 

In view of Theorem \ref{IkawaTasaki} it therefore remains to investigate the maximal Lie subgroups of $G$. The connected maximal Lie subgroups of compact simple Lie groups are well known from classical theory. Due to connectedness we can equivalently consider maximal subalgebras of compact simple Lie algebras. In Table \ref{maxsubalgebra} we list the maximal subalgebras of minimal codimension in compact simple Lie algebras (see e.g.\ \cite{Ma}). 

\begin{table}[h]
\caption{Maximal subalgebras $\ch$ of minimal codimension $d(\cg)$ in compact simple Lie algebras $\cg$} 
\label{maxsubalgebra} 
{\footnotesize\begin{tabular}{ | p{2cm}  p{2cm}  p{2cm}   |}
\hline \rule{0pt}{4mm}
\hspace{-1mm}$\cg$ & $\ch$ & $d(\cg)$  \\[1mm]
\hline \rule{0pt}{4mm}
\hspace{-2mm} 
$\cs\cu_{r+1}$ & $\cs\cu_r \oplus \RR$ & $2r$ \\
$\cs\co_{2r+1}$ & $\cs\co_{2r}$ & $2r$  \\
$\cs\cp_r$ & $\cs\cp_{r-1} \oplus \cs\cp_1$ & $4r-4$  \\
$\cs\co_{2r}$ & $\cs\co_{2r-1}$ & $2r-1$  \\
$\ce_6$ & $\cf_4$ & $26$  \\
$\ce_7$ & $\ce_6 \oplus \RR$ & $54$  \\
$\ce_8$ & $\ce_7 \oplus \cs\cp_1$ & $112$  \\
$\cf_4$ & $\cs\co_9$ & $16$  \\
$\cg_2$ & $\cs\cu_3$ & $6$ \\[1mm]
\hline
\end{tabular}}
\end{table}

We can now finish the proof of Theorem \ref{main}. From Tables \ref{reflindex} and \ref{maxsubalgebra} we get $i_r(G) \leq d(\cg)$. Theorem \ref{IkawaTasaki} then implies $i(G) = i_r(G)$. Using Table \ref{reflindex} we obtain the column for $i(G)$ in Table \ref{Liegroup}.

To find all $\Sigma$ in $G$ with $\codim(\Sigma) = i(G)$ we first note that $i(G) < d(\cg)$ if and only if $G \in \{SU_2,SU_3\}$. In this case $\Sigma$ must be a Cartan embedding and hence a reflective submanifold. From Table \ref{reflindex} we obtain that $\Sigma = SU_2/S(U_1U_1)$ if $G = SU_2$ and $\Sigma = SU_3/SO_3$ if $G = SU_3$. Now assume that $i(G) = d(\cg)$. Then $\Sigma$ is either a Cartan embedding (and then $\Sigma$ is as in Table \ref{reflindex}) or a maximal connected subgroup $H$ of $G$ for which $\ch$ has minimal codimension $d(\cg)$ (and then $\ch$ is as in Table \ref{maxsubalgebra}). By inspection we see that such $H$ is reflective unless $G = G_2$, in which case we get the non-reflective totally geodesic submanifold $SU_3$ of $G_2$ satisfying $\codim(SU_3) = 6 = i(G_2)$. This finishes the proof of Theorem \ref{main}. 

Regarding our conjecture $i(M) = i_r(M)$ if and only if $M \neq G_2^2/SO_4$, we list in Table \ref{summary} the irreducible Riemannian symmetric spaces of noncompact type for which the conjecture remains open.

\begin{table}[h]
\caption{The reflective index $i_r(M)$ for irreducible Riemannian symmetric spaces $M$ of noncompact type for which the conjecture $i(M) = i_r(M)$ is still open and reflective submanifolds $\Sigma$ of $M$ with $\codim(\Sigma) = i_r(M)$} 
\label{summary} 
{\footnotesize\begin{tabular}{ | p{2.9cm}  p{3.7cm}  p{1.5cm}  p{0.8cm}  p{2cm}  |}
\hline \rule{0pt}{4mm}
\hspace{-1mm}$M$ & $\Sigma$ & $\dim M$ & $i_r(M)$ & Comments  \\[1mm]
\hline \rule{0pt}{4mm}
\hspace{-2mm} 
$SU^*_{2r+2}/Sp_{r+1}$ & $\RR \times SU^*_{2r}/Sp_r$ & $r(2r+3)$ & $4r$ & $r \geq 3$ \\[1mm]
\hline \rule{0pt}{4mm}
\hspace{-2mm} 
$Sp_r(\RR)/U_r$ & $\RR H^2 \times Sp_{r-1}(\RR)/U_{r-1}$ & $r(r+1)$ & $2r-2$ & $r \geq 6$\\
$SO^*_{4r}/U_{2r}$ & $SO^*_{4r-2}/U_{2r-1}$ & $2r(2r-1)$ & $4r-2$ &  $r \geq 3$ \\
$Sp_{r,r}/Sp_rSp_r$ & $Sp_{r-1,r}/Sp_{r-1}Sp_r$ &  $4r^2$ & $4r$ & $r \geq 3$ \\
$E_7^{-25}/E_6U_1$ & $E_6^{-14}/Spin_{10}U_1$ & $54$ & $22$ &  \\[1mm]
\hline \rule{0pt}{4mm}
\hspace{-2mm} 
$Sp_{r,r+k}/Sp_rSp_{r+k}$ & $Sp_{r,r+k-1}/Sp_rSp_{r+k-1}$ & $4r(r+k)$ & $4r$ &  $r \geq 3, k \geq 1$, $r > k+1$ \\
$SO^*_{4r+2}/U_{2r+1}$ &$SO^*_{4r}/U_{2r}$ & $2r(2r+1)$ & $4r$ & $r \geq 3$ \\[1mm]
\hline \rule{0pt}{4mm}
\hspace{-2mm}  
$E_6^6/Sp_4$ & $F_4^4/Sp_3Sp_1$ &  $42$ & $14$ & \\[1mm]
\hline \rule{0pt}{4mm}
\hspace{-2mm} 
$E_7^7/SU_8$ & $\RR \times E^6_6/Sp_4$ & $70$ & $27$ & \\[1mm]
\hline \rule{0pt}{4mm}
\hspace{-2mm} 
$E_8^8/SO_{16}$ & $\RR H^2 \times E_7^7/SU_8$ & $128$ & $56$ & \\[1mm]
\hline \rule{0pt}{4mm}
\hspace{-2mm} 
$E_6^2/SU_6Sp_1$ & $F_4^4/Sp_3Sp_1$ & $40$ & $12$ &  \\
$E_7^{-5}/SO_{12}Sp_1$ & $E_6^2/SU_6Sp_1$ & $64$ & $24$ & \\
$E_8^{-24}/E_7Sp_1$ & $E_7^{-5}/SO_{12}Sp_1$ & $112$ & $48$ & \\[1mm]
\hline
\end{tabular}}
\end{table}

\end {document}